\documentclass[aos,MSNbibl,citesort,dvips]{arximspdf}

\doi{10.1214/11-AOS959}
\volume{40}
\issue{1}
\pubyear{2012}
\firstpage{385}
\lastpage{411}

\makeatletter

\renewcommand{\epsilon}{\varepsilon}
\newcommand{\vvert}{\Vert}

\newtheorem{theorem}{Theorem}[section]

\newtheorem{lemma}[theorem]{Lemma}

\newproclaim{remark}{Remark}[section]

\newcommand{\E}{{\mathbb E}}
\newcommand{\R}{{\mathbb R}}
\renewcommand{\P}{{\mathbb P}}
\newcommand{\Q}{{\mathbb Q}}
\newcommand{\K}{{\mathcal{K}}}
\newcommand{\Ps}{{\mathcal{P}}}
\newcommand{\F}{{\mathcal F}}
\newcommand{\scap}{{\mathcal{S}}}
\newcommand{\ham}{{\Upsilon}}
\newcommand{\dok}{{\breve{K}}}
\newcommand{\rok}{{\breve{R}}}
\newcommand{\euc}{{\mathfrak{d}}}
\newcommand{\unisp}{{\nu_{\mathrm{unif}}}}

\newcommand{\argmin}{\mathop{\arg\min}}

\makeatother

\begin{document}
\begin{frontmatter}

\title{Optimal rates of convergence for convex set estimation from support functions}
\runtitle{Convex set estimation from support functions}

\begin{aug}
\author[A]{\fnms{Adityanand} \snm{Guntuboyina}\corref{}\ead[label=e1]{aditya@stat.berkeley.edu}\ead[label=u1,url]{http://www.stat.berkeley.edu/\textasciitilde aditya}}
\runauthor{A. Guntuboyina}
\affiliation{University of Pennsylvania}
\address[A]{Department of Statistics\\
University of California, Berkeley\\
423 Evans Hall\\
Berkeley, California 94720\\
USA\\
\printead{e1}\\
\printead{u1}} 
\end{aug}

\received{\smonth{8} \syear{2011}}
\revised{\smonth{11} \syear{2011}}

%
\begin{abstract}
We present a~minimax optimal solution to the problem of estimating a~compact, convex set from finitely many noisy measurements of its
support function. The solution is based on appropriate regularizations
of the least squares estimator. Both fixed and random designs are
considered.
\end{abstract}

%

\end{frontmatter}

\section{Introduction}
We study the nonparametric estimation problem of estimating a~compact, convex set in Euclidean space from noisy support function
measurements. The support function $h_K$ of a~compact, convex subset
$K$ of $\R^d$ ($d \geq2$) is defined for $u$ in the unit sphere,
$S^{d-1} := \{x \in\R^d \dvtx x_1^2 + \cdots+ x_d^2 = 1 \}$ by
\[
h_K(u) := \sup_{x \in K} ( x \cdot u ) \qquad\mbox{where $x \cdot u :=
x_1u_1 +\cdots+ x_du_d$}.
\]
The function $h_K$ is called the support function of $K$ because it
provides information on support hyperplanes and halfspaces of
$K$. Indeed, every support halfspace of $K$ is of the form $\{x \dvtx x
\cdot u \leq h_K(u)\}$ for some $u \in S^{d-1}$ and since $K$ equals
the intersection of all its support halfspaces, the function $h_K$
uniquely determines $K$. For a~proof of this and other elementary
properties of the support function, see Schneider~\cite{Schneider},
Section 1.7, or Rockafellar~\cite{Rockafellar70book},
Section 13.

We consider the problem of estimating an unknown compact, convex
set~$K$ from observations $(u_1, Y_1),\ldots, (u_n, Y_n)$ drawn according
to the model
%
\begin{equation}\label{arunaroy}
Y_i = h_K(u_i) + \xi_i \qquad\mbox{for $i = 1,\ldots, n$},
\end{equation}
where $u_1,\ldots, u_n$ are unit\vspace*{1pt} vectors and $\xi_1,\ldots, \xi_n$ are
independent normally distributed random variables with mean zero and
variance $\sigma^2$. We work with both fixed and random-design
settings for $u_1,\ldots, u_n$.

Direct motivation for a~theoretical study of this problem comes from
applications. For example, Prince and Willsky~\cite{PrinceWillskyIEEE}, who
were the first to propose the regression model~(\ref{arunaroy}) for
this problem, were motivated by an application to Computed
Tomography. Lele, Kulkarni and Willsky
\cite{LeleKulkarniWillsky} showed how solutions to this
problem can be applied to target reconstruction from resolved
laser-radar measurements in the presence of registration
errors. Gregor and Rannou~\cite{GregorRannou} considered applications
to Projection
Magnetic Resonance Imaging. Another application domain where this
problem might plausibly arise is robotic tactile sensing as has been
suggested by Prince and Willsky~\cite{PrinceWillskyIEEE}. Under an
observation model
that is different from the one considered
here, Goldenshluger and Zeevi~\cite{GoldenshlugerZeevi2006} studied
the estimation of the
support function of the convex support of an unknown intensity
function in the context of image analysis.

Additional reasons for analyzing this estimation problem arise from
the fact that it has a~similar flavor to other well-studied
regression problems:
\begin{longlist}[2.]
\item[1.] It is essentially a~nonparametric function estimation problem
where the true function is assumed to be the support function of a~compact, convex set, that is, there is an implicit convexity-based
constraint on the true regression function. Regression and density
estimation problems with explicit such constraints, for example,
log-concave density estimation and convex regression, have received
much attention. Some examples of work in this general area include
Balabdaoui and Wellner~\cite{BalabdaouiWellner2007}, Balabdaoui,
Rufibach and Wellner~\cite{BalabdaouiRufibachWellner2009}, Cule,
Samworth and Stewart~\cite{CuleSamworthStewart}, D\"umbgen and Rufibach
\cite{DumbgenRufibach2009}, Groeneboom, Jongbloed and Wellner
\cite{GroeneboomJongbloedWellner2001a,GroeneboomJongbloedWellner2001b},
Mammen~\cite{mammen91AS}, Seijo and Sen~\cite{SeijoSen2011}.

\item[2.] Our model $Y_i = \max_{x \in K} (x \cdot u_i) + \xi_i$ can
also be viewed as a~variant of the usual linear regression model
where the dependent variable is modeled as the maximum of linear
combinations of the explanatory variables over a~set of parameter
values and where the interest lies in estimating the convex hull of
the set of parameters. While we do not know if this maximum
regression model has been used outside the context of convex set
estimation, the idea of combining linear functions of independent
variables into nonlinear algorithmic prediction models for the
response variable is familiar (as in neural networks).
\end{longlist}

Let us now briefly describe the previous work on this problem. The
least squares estimator has been the most commonly used. It
is defined as
%
\begin{equation}\label{ls}
\hat{K}_{\mathrm{ls}} := \argmin_{L} \sum_{i=1}^n \bigl(Y_i - h_L(u_i)
\bigr)^2,
\end{equation}
where the minimum is taken over all compact, convex subsets $L$. The
minimizer here is not unique and one can always take it to be a~polytope (convex set with finitely many corners; more carefully
defined in the next section). This estimator, for\vadjust{\goodbreak}
$d=2$, was first proposed by Prince and Willsky~\cite{PrinceWillskyIEEE}, who assumed
that $u_1,\ldots, u_n$ are evenly spaced on the unit circle and that
the error variables $\xi_1,\ldots, \xi_n$ are normal with zero
mean. They also proposed an algorithm for computing it based on
quadratic programming. Lele, Kulkarni and Willsky
\cite{LeleKulkarniWillsky} extended this
algorithm to include the case of nonevenly spaced two-dimensional
$u_1,\ldots, u_n$ as well. Recently, Gardner and Kiderlen~\cite{GardnerKiderlen2009}
proposed an algorithm for computing a~minimizer of the least squares
criterion for every dimension $d \geq2$ and every sequence $u_1,\ldots, u_n$.

In addition to the least squares estimator, Prince and Willsky~\cite{PrinceWillskyIEEE}
and
Lele, Kulkarni and Willsky~\cite{LeleKulkarniWillsky} also proposed
estimators (in the case
\mbox{$d=2$}) designed to take advantage of certain forms of prior knowledge,
when available, about the true compact, convex set. These estimators
are all based on a~least squares minimization.

Fisher, Hall, Turlach and Watson
\cite{FisherHallTurlachWatson} proposed estimators for $d=2$ that are
not based on the least squares criterion. They made smoothness
assumptions on the true support function $h_K$ (viewed as a~function
on the unit circle or on the interval $(-\pi, \pi]$) and estimated it
using periodic versions of standard nonparametric regression
techniques such as local regression, kernel smoothing and
splines. They suggested a~way to convert the estimator of $h_K$ into
an estimator for $K$ using a~formula, which works for smooth~$h_K$,
for the boundary of $K$ in terms of $h_K$. Hall and Turlach
\cite{HallTurlach} added a~corner-finding technique to the method
of Fisher et al.~\cite{FisherHallTurlachWatson} to estimate
two-dimensional convex
sets with certain types of corners.

There are relatively fewer theoretical results in the
literature. Fisher et al.~\cite{FisherHallTurlachWatson}, Theorem 4.1,
stated a~theorem without
proof which appears to imply consistency and certain rates of
convergence for their estimator under certain smoothness assumptions
on the true support function.
Gardner, Kiderlen and Milanfar~\cite{GardnerKiderlenMilanfar} proved
consistency of the
least squares estimator and also derived rates of convergence. They
worked with the following assumptions:
\begin{longlist}[3.]
\item[1.] $u_1, u_2, \ldots$ are deterministic satisfying
\[
\max_{u \in S^{d-1}} \min_{1 \leq i \leq n} \euc_d(u, u_i) =
O\bigl(n^{-1/(d-1)}\bigr) \qquad\mbox{as $n \rightarrow\infty$},
\]
where $\euc_d$ denotes the usual Euclidean distance on $\R^d$
[see~(\ref{l12})].
\item[2.] $\xi_1, \xi_2, \ldots$ are independently distributed according to
the normal distribution with mean zero and variance $\sigma^2$.
\item[3.] $K$ is contained\vspace*{-2pt} in the ball of radius $\Gamma$ centered at the
origin with $\Gamma\geq\sigma^{15/2}$.
\end{longlist}
For the loss function $\ell_f^2$ defined in~(\ref{fdlf})
below, Gardner et al.~\cite{GardnerKiderlenMilanfar}, Corollary~5.7,
showed that
$\ell_f^2(K, \hat{K}_{\mathrm{ls}}) = O_{d, \sigma, \Gamma}(\beta_n)$ as~$n$
approaches $\infty$ almost surely, where
%
\begin{equation}\label{gardls}
\beta_n := \cases{
n^{-4/(d+3)}, &\quad when $d =
2,3,4$,\vspace*{1pt}\cr
n^{-1/2} (\log n )^2, &\quad when
$d = 5$,\vspace*{2pt}\cr
n^{-2/(d-1)}, &\quad
when $d \geq6$.}
\end{equation}
Here $O_{d,\sigma, \Gamma}$ is the usual big-O notation where the
constant involved depends on $d, \sigma$ and
$\Gamma$. Gardner et al.~\cite{GardnerKiderlenMilanfar}, Corollary
5.7, provided
explicit expressions for the dependence of the constant with respect
to $\sigma$ and $\Gamma$ (but not~$d$) which we have suppressed here
because our interest only lies in the dependence on~$n$. Also
see Gardner et al.~\cite{GardnerKiderlenMilanfar}, Lemma 3.2, for implications
of Assumption~1 on the unit vector sequence $\{u_i \}$.

Such a~strange dependence of the rates of convergence of the least
squares estimation on dimension has also been observed in other
situations (see, e.g.,
Birg\'e and Massart~\cite{BirgeMassartMCE},
van de Geer~\cite{VandegeerBook}, Seregin and Wellner~\cite{SereginWellner}).

Our results in this paper, described below, imply that
the rates~(\ref{gardls}) proved
by Gardner et al.~\cite{GardnerKiderlenMilanfar} for the least squares
estimator are
optimal when $d \leq4$ and suboptimal when $d \geq5$. We show how
estimators can be constructed that converge at the rate $n^{-4/(d+3)}$
for all dimensions $d \geq2$. Our estimators are based on
regularizing $\hat{K}_{\mathrm{ls}}$ by minimizing the least squares criterion
on certain well-chosen subsets of the parameter space. In contrast
to Gardner et al.~\cite{GardnerKiderlenMilanfar}, we took the more customary
approach in nonparametric statistics by proving rates for the expected
loss or risk instead of almost sure bounds for the loss. An advantage
is that this results in bounds for a~finite (though large)~$n$ thereby
circumventing the need to let~$n$ approach infinity.

We establish an optimal minimax theory for this problem in both
fixed-design and random-design frameworks.

In the fixed-design framework, we assume that $u_1,\ldots, u_n$ are
deterministic. We define the minimax risk in this setting as (the
subscript $f$ below stands for fixed-design):
%
\begin{equation}\label{fdmr}
R_f(n) = R_f(n; \sigma, \Gamma) := \inf_{\hat{K}} \sup_{K \in
\K^d(\Gamma)} \E_K \ell_f^2(K, \hat{K})
\end{equation}
with
%
\begin{equation}\label{fdlf}
\ell^2_f(K, K') := \frac{1}{n} \sum_{i=1}^n \bigl(h_K(u_i) -
h_{K'}(u_i) \bigr)^2,
\end{equation}
where $\K^d(\Gamma)$ denotes the set of all compact, convex sets
contained in the ball of radius $\Gamma$ centered at the origin and
$\E_K$ denotes expectation taken when the true compact, convex set
equals $K$. We assume that $\sigma$ and $\Gamma$ are known. The
infimum in the definition of $R_f(n)$ is over all possible estimators~$\hat{K}$ where estimators are defined to be functions of $(u_1, Y_1),\ldots, (u_n, Y_n)$ as well as of $\sigma$ and $\Gamma$ taking values
in the space of all compact, convex sets.

For every deterministic set of unit vectors $u_1,\ldots, u_n$, we show
that $R_f(n)$ is bounded from above by a~constant (which is independent
of~$n$) multiple of $n^{-4/(d+3)}$. Under a~specific assumption on
$u_1,\ldots, u_n$, we also prove that a~constant multiple of
$n^{-4/(d+3)}$ is a~lower bound for $R_f(n)$. The upper bound is proved
by considering least squares estimators\vspace*{1pt} on appropriate
subsets of the set of all compact, convex sets in $\R^d$. The lower
bound is proved by the application of Assouad's lemma to a~special
finite collection of convex sets.\vadjust{\goodbreak}

We also study the random-design setting where we assume that $u_1,\ldots, u_n$ are independently distributed according to a~fixed
probability measure, $\nu$, on $S^{d-1}$ and that the errors $\xi_1,\ldots, \xi_n$ are independent of $u_1,\ldots, u_n$. Here we define the
minimax risk as (the subscript $r$ below stands for random design):
%
\begin{equation}\label{rdmr}
R_r(n) = R_n(n; \sigma, \Gamma) := \inf_{\hat{K}} \sup_{K \in
\K^d(\Gamma)} \E_K\ell_r^2(K, \hat{K})
\end{equation}
with
%
\begin{equation}\label{rdlf}
\ell_r^2(K, K') := \int_{S^{d-1}} \bigl(h_K(u) - h_{K'}(u) \bigr)^2
\,d\nu(u).
\end{equation}
For every probability measure $\nu$ on $S^{d-1}$, we show that
$R_r(n)$ is bounded from above by a~constant (which is independent on
$n$) multiple of $n^{-4/(d+3)}$. The proof techniques here are similar
to the fixed-design setting. When~$\nu$ equals $\unisp$, the uniform
probability measure on $S^{d-1}$, we prove that a~constant multiple of
$n^{-4/(d+3)}$ is also a~lower bound for $R_r(n)$. We use a~different
lower bound proof here from the one used in the fixed-design setting.

We would like to remark here that the rate $n^{-4/(d+3)}$ is quite
natural in connection to the minimax estimation of smooth
functions. Indeed, the unit sphere has dimension $d-1$, and the class
of smooth functions on a~space of dimension $d-1$ with smoothness
$\gamma$ allows the minimax rate $n^{-2\gamma/(2\gamma+ d - 1)}$. Our
problem here has a~convexity constraint and convexity is associated,
in a~broad sense, with the smoothness $\gamma= 2$, which explains the
rate $n^{-4/(d+3)}$.

After setting up the necessary notation in the next section, we prove
the fixed-design bounds in Section~\ref{fd} and the random-design
bounds in Section~\ref{rd}. Some auxiliary results that are needed for
the proofs of the main theorems are collected in the three
\hyperref[app]{Appendices}.

\section{Notation}\label{stanot}
This section will fix notation and introduce some standard notions
that are used in the paper.

$\P_K$ denotes the probability distribution of the observations when
the true compact, convex set equals $K$. In other words, $\P_K$ is the
joint distribution of $(Y_1,\ldots, Y_n)$ in the fixed-design setting
and the joint distribution of $(u_1, Y_1),\ldots, (u_n, Y_n)$ in the
random-design setting. We use the same notation in both cases as the
setting will be clear from the context. Expectation under $\P_K$ is
denoted by $\E_K$.

For a~real-valued function $f$ on $S^{d-1} \times\cdots\times
S^{d-1}$, let $\E_{\nu} f(u_1,\ldots, u_n)$ denote expectation taken
under the assumption that $u_1,\ldots, u_n$ are independently
distributed according to $\nu$.

We denote the usual Euclidean distance on $\R^d$ by $\euc_d$, that
is, for
$x = (x_1,\ldots, x_d)$ and $y = (y_1,\ldots, y_d)$ in $\R^d$,
%
\begin{equation}\label{l12}
\euc_d(x, y) := \Biggl(\sum_{i=1}^d (x_i - y_i)^2 \Biggr)^{1/2}.
\end{equation}
The closed ball in $\R^d$ of radius $\alpha> 0$ centered at $a~\in
\R^d$ is denoted by $B_d(a, \alpha)$, that is, $B_d(a, \alpha) := \{
x \in
\R^d \dvtx \euc_d(x, a) \leq\alpha\}$.

The uniform probability measure on $S^{d-1}$ is denoted by $\unisp$.

The convex hull of a~subset $S$ of $\R^d$, defined as the intersection
of all convex sets containing $S$, is denoted by
$\operatorname{conv}(S)$. Convex hulls of finite subsets of $\R^d$
are called
polytopes. The set of all polytopes in $\R^d$ with at most $m$ extreme
points (corners) is denoted by $\Ps_m$, that is,
\[
\Ps_m := \{\operatorname{conv}(S) \dvtx S \subseteq\R^d \mbox{ with
cardinality at most } m \}.
\]
For $\Gamma> 0$, let $\Ps_m(\Gamma)$ denote the set of all polytopes in
$\Ps_m$ that are contained in $B_d(0, \Gamma)$:
%
\begin{equation}\label{polgam}
\Ps_m(\Gamma) := \{K \in\Ps_m \dvtx K \subseteq B_d(0,
\Gamma)\}.
\end{equation}

For two compact, convex subsets $K$ and $K'$ of $\R^d$, the Hausdorff
distance between them is defined as
%
\begin{equation}\label{haus1}
\ell_H(K, K') := {\sup_{u \in S^{d-1}}} |h_K(u) - h_{K'}(u)|.
\end{equation}
It is apparent that both $\ell_f^2(K, K')$ and $\ell_r^2(K, K')$ are
less than or equal to $\ell_H^2(K, K')$. The Hausdorff distance has
the following alternative expression:
%
\begin{equation}
\label{althaus}
\ell_H(K, K') = \max\Bigl(\sup_{x \in K} \inf_{y \in K'}
\euc_d(x, y), \sup_{x \in K'} \inf_{y \in K}
\euc_d(x, y) \Bigr).
\end{equation}
A~simple proof of the equivalence of~(\ref{haus1}) and~(\ref{althaus})
can be found in Schneider~\cite{Schneider}, Theorem 1.8.11.

The standard notions of packing and covering numbers will be frequently
used and we have collected their definitions below for the convenience
of the reader. Let $\Theta$ be an arbitrary set and
let $\rho$ be a~nonnegative function on $\Theta\times\Theta$ ($\rho$
does not necessarily have to be a~metric).
\begin{longlist}[2.]
\item[1.] \textit{Packing}: By an $\eta$-packing subset of
$(\Theta, \rho)$, we mean a~subset $S \subseteq\Theta$ for which
$\rho(\theta, \theta') \geq\eta$ for all $\theta, \theta' \in S$
with $\theta\neq
\theta'$. By a~maximal $\eta$-packing subset, we mean an
$\eta$-packing subset that is not a~proper subset of any other
$\eta$-packing subset. The packing number $N(\Theta, \eta; \rho)$ is
defined as the maximum of the cardinalities of all $\eta$-packing
subsets of $\Theta$.

\item[2.] \textit{Covering}: By an $\epsilon$-covering subset of $(\Theta,
\rho)$, we mean a~subset
$S \subseteq\Theta$ such that $\min_{s \in S} \rho(t, s) \leq
\epsilon$ for every $t \in\Theta$. The $\epsilon$-covering number
$M(\Theta, \epsilon; \rho)$ is defined as the minimum of the
cardinalities of all $\epsilon$-covering subsets of $\Theta$.
\end{longlist}

The cardinality of a~finite set $F$ is denoted by $|F|$.

We use the following notions of distance between probability
measures $P$ and $Q$ having densities $p$ and $q$ with respect to a~common measure $\mu$:
\begin{longlist}[3.]
\item[1.] Total variation distance: $\|P-Q\|_{\mathrm{TV}} := \int(|p - q|/2)
\,d\mu$.
\item[2.] Kullback--Leibler divergence: $D(P\vvert Q) := \int p \log(p/q) \,d\mu$
if $P$ is absolutely continuous with respect to $Q$ and $\infty$
otherwise.\vadjust{\goodbreak}
\item[3.] Chi-squared divergence: $\chi^2(P\vvert Q) := \int(p^2/q) \,d\mu- 1$
if $P$ is absolutely continuous with respect to $Q$ and $\infty$
otherwise. We also write $\chi(P, Q) := ( \chi^2(P\vvert Q)
)^{1/2}$.
\end{longlist}
Pinsker's inequality states
%
\begin{equation}\label{pinq}
D(P\vvert Q) \geq2 \|P - Q\|^2_{\mathrm{TV}}
\qquad\mbox{for all probability measures $P$
and $Q$}.
\end{equation}

We use the symbols $c, C$, etc. to denote positive constants depending
only on the dimension~$d$. Their value may change with every
occurrence.

\section{Fixed design setting}\label{fd}
In this section, we assume that $u_1,\ldots, u_n$ are
deterministic. The errors $\xi_1,\ldots, \xi_n$ are assumed to be
independently distributed according to the normal distribution with
mean zero and variance~$\sigma^2$. We consider the loss function
$\ell_f^2$ [defined in~(\ref{fdlf})] and prove upper and lower bounds
for the corresponding minimax risk $R_f(n)$ over $\K^d(\Gamma)$ [see
the definition~(\ref{fdmr})].

\subsection{Upper bound for $R_f(n)$}
The following result shows that $R_f(n)$ is at most $n^{-4/(d+3)}$ up
to a~multiplicative constant that is independent of~$n$. We make no
assumptions on the deterministic design unit vectors $u_1,\ldots, u_n$
and they are completely arbitrary. It should be noted that the loss
function~$\ell_f^2$ is naturally associated with this fixed design
setup enabling the following theorem to hold with no assumptions
whatsoever on $u_1,\ldots, u_n$. On the other hand, such assumptions
would be unavoidable if one is interested in proving risk bounds for
other loss functions under the fixed design setting. A~similar\vspace*{1pt} remark also applies to Theorem~\ref{ubrthm} where
the natural loss function is $\ell_r^2$.
%
\begin{theorem}\label{ubfthm}
There exist positive constants $c$ and $C$ depending only on the
dimension~$d$ such that
%
\begin{equation}\label{ubfthmeq}
R_f(n) \leq c \sigma^{8/(d+3)} \Gamma^{2(d-1)/(d+3)} n^{-4/(d+3)}
\qquad\mbox{if $n \geq C(\sigma/\Gamma)^2$}.
\end{equation}
\end{theorem}
\begin{pf}
For each finite subset $F$ of $\K^d(\Gamma)$, let us define the least
squares estimator $\hat{K}_F$ by
%
\begin{equation}\label{fls}
\hat{K}_F := \argmin_{L \in F} \sum_{i=1}^n \bigl(Y_i - h_L(u_i)
\bigr)^2.
\end{equation}
We show that, if $F$ is chosen appropriately, then $\E_K \ell_f^2(K,
\hat{K}_F)$ is bounded from above by the right-hand side
of~(\ref{ubfthmeq}) for every $K \in\K^d(\Gamma)$.

Fix $K \in\K^d(\Gamma)$. We start with the following trivial
inequality which holds for every nonnegative function $G$ on $F$ and
every $\alpha> 0$:
\[
G(\hat{K}_F) \leq\sum_{L \in F} G(L) \exp\biggl(\alpha\sum_i
\bigl(Y_i - h_{\hat{K}_F}(u_i)\bigr)^2 - \alpha\sum_i
\bigl(Y_i - h_L(u_i)\bigr)^2 \biggr),
\]
the reason being that the term for $L = \hat{K}_F$ in the sum on the
right-hand side equals $G(\hat{K}_F)$.

Because $\hat{K}_F$ is the least squares estimator over $F$, we can
replace it in the right-hand side above by an arbitrary $L' \in
F$. Taking expectation on both sides of the resulting inequality, we
obtain
%
\begin{equation}\label{repineq}
\E_K G(\hat{K}_F) \leq\sum_{L \in F} G(L) \E_K e^{\alpha
\sum_i (Y_i - h_{L'}(u_i))^2 - \alpha\sum_i (Y_i - h_L(u_i))^2}
\end{equation}
for every $L' \in F$. The expectation term in the right-hand side
equals
\[
\exp\bigl(-\alpha n \ell_f^2(K, L) + \alpha n \ell_f^2(K, L')
+ 2 \alpha^2 \sigma^2 n \ell_f^2(L, L') \bigr).
\]
We then use the elementary inequality $\ell_f^2(L, L') \leq2
\ell_f^2(K, L) + 2 \ell_f^2(K, L')$ and the fact that $\ell_f^2(K, L')
\leq\ell_H^2(K, L')$ to obtain the following upper bound for $\E_K
G(\hat{K}_F)$:
\[
\min_{L' \in F} \sum_{L \in F} G(L) \exp
\bigl((-\alpha+4\alpha^2\sigma^2) n \ell_f^2(K, L) +
(\alpha+4\alpha^2\sigma^2) n \ell_H^2(K, L') \bigr).
\]
The choices
\[
G(L) = \exp\biggl(\frac{n\ell_f^2(K, L)}{16 \sigma^2} \biggr)
\quad\mbox{and}\quad \alpha= \frac{1}{8 \sigma^2}
\]
lead to the following risk bound:
%
\begin{equation}\label{expb}
\E_K \exp\biggl(\frac{n \ell_f^2(K, \hat{K}_F)}{16 \sigma^2}
\biggr)
\leq|F| \exp\biggl(\frac{3n}{16 \sigma^2} \min_{L' \in F}
\ell_H^2(K, L') \biggr),
\end{equation}
where $|F|$ denotes the cardinality of $F$. Using Jensen's inequality
on the left-hand side and taking logarithms on both sides, we deduce
\[
\E_K \ell_f^2(K, \hat{K}_F) \leq{\frac{16 \sigma^2}{n} \log}|F| + 3
\min_{L' \in F} \ell_H^2(K, L').
\]
Since $K \in\K^d(\Gamma)$ was arbitrary, we get
%
\begin{equation}\label{bhupali}
R_f(n) \leq{\frac{16 \sigma^2}{n} \log}|F| + 3
\sup_{K \in\K^d(\Gamma)} \min_{L' \in F} \ell_H^2(K, L').
\end{equation}
We now use a~classical result on the covering numbers of
$(\K^d(\Gamma), \ell_H)$ due to Bronshtein~\cite{Bronshtein76}, Theorem 3 and Remark
1, which states that there exist postive constants $c$
and $\epsilon_0$ which depend on~$d$ alone such that for every
$\epsilon\leq\Gamma\epsilon_0$, there exists a~finite subset $F
\subseteq\K^d(\Gamma)$ satisfying
%
\begin{equation}\label{bron1}
{\log}|F| \leq c \biggl(\frac{\Gamma}{\epsilon} \biggr)^{(d-1)/2}
\quad\mbox{and}\quad \sup_{K \in K^d(\Gamma)} \min_{L' \in F} \ell_H^2(K,
L') \leq\epsilon^2.
\end{equation}
Combining~(\ref{bhupali}) and~(\ref{bron1}), we get
\[
R_f(n) \leq\biggl[16 c \frac{\sigma^2}{n}
\biggl(\frac{\Gamma}{\epsilon} \biggr)^{(d-1)/2} + 3 \epsilon^2
\biggr] \qquad\mbox{for every $0 < \epsilon\leq\Gamma\epsilon_0$.}
\]
If we now choose $\epsilon:= \sigma^{4/(d+3)} \Gamma^{(d-1)/(d+3)}
n^{-2/(d+3)}$, then $\epsilon\leq\Gamma\epsilon_0$ provided $n \geq
C(\sigma/\Gamma)^2$ for a~large enough constant $C$ depending only on
$d$ and the required inequality~(\ref{ubfthmeq}) follows.
\end{pf}

\subsection{Lower bound for $R_f(n)$}
We show that $n^{-4/(d+3)}$ is also a~lower bound for $R_f(n)$ up to a~multiplicative constant that is independent of~$n$. We make the
assumption that the fixed unit vectors $u_1,\ldots, u_n$ form a~maximal $\epsilon$-packing subset (under the Euclidean metric
$\euc_d$) of $S^{d-1}$ for some $\epsilon\in(0, 1]$. The definition
of a~maximal packing set was given in Section~\ref{stanot}. Note that
it is impossible to prove the lower bound $n^{-4/(d+3)}$ for $R_f(n)$
without any assumptions on $u_1,\ldots, u_n$. For example, if $u_1 =
\cdots= u_n$, then $R_f(n)$ is of the order $1/n$.

A~standard argument [sketched in the \hyperref[app]{Appendix}; see
inequality~(\ref{pent})] shows that our assumption on $u_1,\ldots,
u_n$ implies that
%
\begin{equation}\label{maxpack}
c \epsilon^{1-d} \leq n \leq C \epsilon^{1-d}
\end{equation}
for two constants $c$ and $C$ depending only on~$d$. The following is
the main theorem of this section.
%
\begin{theorem}\label{assf}
Suppose $u_1,\ldots, u_n$ form a~maximal $\epsilon$-packing
subset\break
of~$S^{d-1}$ for some $\epsilon\in(0, 1]$. There exist positive
constants $c$ and $C$ depending only on~$d$ such that
%
\begin{equation}\label{assfeq}
R_f(n) \geq c \sigma^{8/(d+3)} \Gamma^{2(d-1)/(d+3)} n^{-4/(d+3)},
\end{equation}
whenever $n \geq C \max( (\sigma/\Gamma)^2,
(\Gamma/\sigma)^{(d-1)/2})$.
\end{theorem}

Our proof is based on the application of Assouad's lemma to an
explicitly constructed finite subset of $\K^d(\Gamma)$. The following
version of Assouad's lemma is taken from van der Vaart
\cite{vaart98book}, page
347. Recall that $\P_K$ denotes the
probability distribution of the observations when the true compact,
convex set is $K$.
%
\begin{lemma}[(Assouad)]
Let $m$ be a~positive integer and suppose that, for each $\tau\in
\{0, 1\}^m$, there is an associated set $K(\tau)$ in
$\K^d(\Gamma)$. Then the following inequality holds:
\[
R_f(n) \geq\frac{m}{8} \min_{\tau\neq\tau'}
\frac{\ell_f^2(K(\tau), K(\tau'))}{\ham(\tau, \tau')} \min
_{\ham
(\tau, \tau') = 1} \bigl(1 - \bigl\|\P_{K(\tau)} - \P_{K(\tau
')}\bigr\|_{\mathrm{TV}} \bigr),
\]
where $\ham(\tau, \tau') := \sum_{i} \{\tau_i \neq\tau'_i\}$.
\end{lemma}
\begin{pf*}{Proof of Theorem~\ref{assf}}
We apply Assouad's lemma to the following construction. For a~fixed
positive $\eta\leq1/8$, we consider unit vectors $v_1,\ldots, v_m$
that form a~maximal $2\sqrt{2\eta}$-packing subset of the unit
sphere (under $\euc_d$) and define
\[
K(\tau) := D_1(\tau_1) \cap\cdots\cap D_m(\tau_m) \qquad\mbox{for $\tau
\in\{0, 1\}^m$},
\]
where
\[
D_i(0) := B_d(0, \Gamma) \cap\{x \dvtx x \cdot v_j \leq
\Gamma(1 - \eta)\} \quad\mbox{and}\quad D_i(1) := B_d(0, \Gamma).
\]
One consequence of the assumption on $v_1,\ldots, v_m$ is that $m \geq
c \eta^{(1-d)/2}$ for a~constant $c$. Another consequence is that the
sets $B_d(0, \Gamma) \cap\{x \cdot v_j > 1 - \eta\}$ are disjoint
which implies that
\[
\ell_f^2(K(\tau), K(\tau')) = \sum_{i \dvtx \tau_i \neq\tau'_i}
\ell_f^2(D_i(0), D_i(1)) = \ham(\tau, \tau') \ell_f^2(D_1(0),
D_1(1))
\]
for every $\tau, \tau' \in\{0, 1\}^m$. In Lemma~\ref{appass} (stated
and proved in Appendix~\ref{appA}), we show that there exist constants $c$ and
$C$ such that
\[
c \Gamma^2 \eta^{(d+3)/2} \leq\ell_f^2(D_1(0), D_1(1)) \leq C
\Gamma^2 \eta^{(d+3)/2}
\]
provided $0 < \eta\leq1/8$ and $\eta\geq C \epsilon^2$. Therefore,
%
\begin{equation}\label{assreq1}\qquad
\min_{\tau\neq\tau'} \frac{\ell_f^2(K(\tau), K(\tau'))}{\ham
(\tau,
\tau')} \geq c \Gamma^2 \eta^{(d+3)/2} \qquad\mbox{if $0 < \eta\leq1/8$
and $\eta\geq C\epsilon^2$}.
\end{equation}
To bound $\|\P_{K(\tau)} - \P_{K(\tau')}\|_{\mathrm{TV}}$, we use Pinsker's
inequality~(\ref{pinq}) because the Kullback--Leibler divergence
$D(P_{K(\tau)}\vvert P_{K(\tau')})$ has a~simple expression in terms of
$\ell_f^2(K(\tau), K(\tau'))$:
\begin{eqnarray*}
\bigl\|\P_{K(\tau)} - \P_{K(\tau')}\bigr\|^2_{\mathrm{TV}} &\leq& \frac{1}{2}
D\bigl(\P_{K(\tau)}\vvert\P_{K(\tau')}\bigr)
= \frac{n}{4\sigma^2}
\ell_f^2(K(\tau), K(\tau'))\\
&=& \frac{n}{4\sigma^2}\Upsilon(\tau,
\tau')
\ell_f^2(D_1(0), D_1(1)).
\end{eqnarray*}
By a~second application of Lemma~\ref{appass}, we obtain
%
\begin{equation}\label{assreq2}
\min_{\Upsilon(\tau, \tau') = 1}\bigl(1 - \bigl\|\P_{K(\tau)} -
\P_{K(\tau')}\bigr\|_{\mathrm{TV}} \bigr) \geq1 - C \frac{\Gamma
\sqrt{n}}{\sigma} \eta^{(d+3)/4},
\end{equation}
if $0 < \eta\leq1/8$ and $\eta\geq C \epsilon^2$. Therefore,
applying Assouad's lemma with the
inequalities~(\ref{assreq1}),~(\ref{assreq2}) and $m \geq
c\eta^{-(d-1)/2}$, we obtain
%
\begin{equation}\label{jaunpuri}
R_f(n) \geq c \Gamma^2 \eta^2 \biggl(1 - C \frac{\Gamma
\sqrt{n}}{\sigma} \eta^{(d+3)/4} \biggr),
\end{equation}
if $0 < \eta\leq1/8$ and $\eta\geq C \epsilon^2$. We now make the
choice
\[
\eta:= \biggl(\frac{1}{2C} \frac{\sigma}{\Gamma\sqrt{n}}
\biggr)^{4/(d+3)}.
\]
Then $\eta\leq1/8$ provided $n \geq C(\sigma/\Gamma)^2$. Also,
from~(\ref{maxpack}), $\epsilon\leq C n^{-1/(d-1)}$ and thus, for
$\eta\geq C \epsilon^2$, it is enough to ensure that $\eta\geq C
n^{-2/(d-1)}$ which, upon simplification, reduces to $n \geq C
(\Gamma/\sigma)^{(d-1)/2}$. Inequality~(\ref{jaunpuri}) with this
choice of $\eta$ then implies~(\ref{assfeq}) which completes the
proof.
\end{pf*}

\subsection{A~more natural estimator}
In Theorem~\ref{ubfthm}, we used the least squares estimator on an
appropriate finite subset of $\K^d(\Gamma)$ to prove the optimal upper
bound for $R_f(n)$. This estimator can be viewed as a~regularized
version of the full least squares estimator $\hat{K}_{\mathrm{ls}} :=
\argmin_{L} \sum_i (Y_i - h_L(u_i))^2$ for which the rates
[see~(\ref{gardls})] proved by Gardner et al.~\cite{GardnerKiderlenMilanfar} are
suboptimal for $d \geq5$.

We remarked [just after~(\ref{ls})] that, for the full least squares
estimator, the set $L$ which minimizes $\sum_i (Y_i - h_L(u_i))^2$ is
not unique. Gardner and Kiderlen~\cite{GardnerKiderlen2009} observed
that a~minimizer can
always be chosen to be a~polytope with at most~$n$ extreme points and
provided an algorithm for computing such a~minimizer. In light of this
observation of Gardner and Kiderlen~\cite{GardnerKiderlen2009}, we
consider\vspace*{1pt} the following
estimator which is a~more intuitive regularization of $\hat{K}_{\mathrm{ls}}$
compared to $\hat{K}_F$:
%
\begin{equation}\label{dokyun}
\dok_m := \argmin_{L \in\Ps_m(\Gamma)} \sum_{i=1}^n \bigl(Y_i -
h_L(u_i) \bigr)^2.
\end{equation}
The set $\Ps_m(\Gamma)$ was defined in~(\ref{polgam}). The best risk
achievable by $\dok_m$ is defined as
\[
\rok_f(n) := \inf_{m \geq1} \sup_{K \in\K^d(\Gamma)} \E_K
\ell_f^2(K, \dok_m).
\]

It is not too hard to see that $\dok_m$ equals the least squares
estimator $\hat{K}_{\mathrm{ls}}$ whenever $m \geq n$. On the other hand, for
$m < n$, they can be quite different.

In this section, we prove the following theorem which shows that
$\rok_f(n)$ is bounded from above by $n^{-4/(d+3)}$ up to a~multiplicative factor that is logarithmic in~$n$. No assumptions on
$u_1,\ldots, u_n$ are necessary.
%
\begin{theorem}\label{lss}
There exist positive constants $c$ and $C$ that depend only on
the dimension~$d$ such that
%
\begin{equation}\label{lsseq}
\rok_f(n) \leq c \sigma^{8/(d+3)}
\Gamma^{2(d-1)/(d+3)}n^{-4/(d+3)} \log(c n \Gamma^2/\sigma^2
) \log(c n),
\end{equation}
if $n \geq C(\sigma/\Gamma)^2$.
\end{theorem}

For the proof of this theorem, we use the following result which is a~special instance of a~result on convergence rates of sieved least
squares estimators from van de Geer~\cite{VandegeerBook}, pages 184
and 185. For a~polytope $P \in\Ps_m$ and $\omega> 0$, let
\[
S_m(P, \omega) := \{L \in\Ps_m \dvtx \ell_f^2(P, L) \leq\omega^2
\}
\]
and let $M(S_m(P, \omega), \epsilon; \ell_f)$ denote the
$\epsilon$-covering number of $S_m(P, \omega)$ under~$\ell_f$.
%
\begin{theorem}\label{vdg}
Fix a~polytope $P \in\Ps_m(\Gamma)$. Suppose $\Psi$ is a~function on
$(0, \infty)$ such that
\[
\Psi(\omega) \geq\int_0^{\omega} \sqrt{\log M(S_m(P, \omega),
\epsilon; \ell_f)} \,d\epsilon\qquad\mbox{for every $\omega> 0$}
\]
and such that $\Psi(\omega)/\omega^2$ is decreasing on $(0,
\infty)$. Then there exists a~universal constant $C$ such that
%
\begin{equation}\label{vdgeq}
\P_K \bigl(\ell_f^2(\dok_m, P) > \delta\bigr) \leq C \sum_{s
\geq
0} \exp\biggl(\frac{-n2^{2s}\delta}{C^2 \sigma^2} \biggr)
\end{equation}
for every $\delta$ satisfying $\delta\geq8 \ell_f^2(K, P)$ and
$\sqrt{n} \delta\geq C \sigma\Psi(\sqrt{\delta})$.
\end{theorem}

The application of this theorem for the proof of~(\ref{lsseq})
requires an upper bound on $M(S_m(P, \omega), \epsilon; \ell_f)$. Such
a~bound is provided in Lemma~\ref{bpw} (stated and proved in Appendix
\ref{appB}).
\begin{pf*}{Proof of Theorem~\ref{lss}}
Fix $m \geq1$ and an arbitrary polytope $P \in\Ps_m(\Gamma)$. In
Lemma~\ref{bpw}, we show that
\[
M(S_m(P, \omega), \epsilon; \ell_f) \leq\biggl(4 +
\frac{2\sqrt{n}\omega}{\epsilon} \biggr)^{b_1 md\log(b_2 m)}
\]
for universal positive constants $b_1$ and $b_2$. This implies that
\begin{eqnarray*}
\int_0^{\omega} \sqrt{\log M(S_m(P, \omega),
\epsilon; \ell_f)} \,d\epsilon&\leq&
\sqrt{b_1md\log(b_2m)} \int_0^{\omega} \sqrt{\log\biggl(4 +
\frac{2 \sqrt{n}\omega}{\epsilon} \biggr)} \,d\epsilon\\
&=& \omega\sqrt{b_1md\log(b_2m)} \int_1^{\infty}
\frac{\sqrt{\log(4+2\sqrt{n}x)}}{x^2}
\,dx \\
&\leq& C \omega\sqrt{md \log(b_2 m) \log\bigl(4+2\sqrt{n}\bigr)}.
\end{eqnarray*}
As a~result, the function $\Psi(\omega)$ appearing in
Theorem~\ref{vdg} can be taken to be
\[
\Psi(\omega) := C \omega\sqrt{md \log(b_2 m) \log\bigl(4+2\sqrt{n}\bigr)}
\qquad\mbox{for every $\omega> 0$}
\]
and then~(\ref{vdgeq}) gives
%
\begin{equation}\label{malkauns}\qquad
\P_K \bigl(\ell_f^2(\dok_m, P) > \delta\bigr) \leq C \sum_{s
\geq
0} \exp\biggl(\frac{-n2^{2s}\delta}{C^2 \sigma^2} \biggr)
\qquad\mbox{whenever $\delta\geq\delta_0$},
\end{equation}
where
\[
\delta_0 := C \biggl(\ell_f^2(K, P) + \frac{\sigma^2}{n}md \log
(b_2m)\log\bigl(4+2\sqrt{n}\bigr) \biggr).
\]
Because $n\delta_0 \geq C \sigma^2$ for a~constant $C$, the sum on the
right-hand side of~(\ref{malkauns}) can be bounded from above by a~constant multiple of the first term, and we deduce
\[
\P_K \bigl(\ell_f^2(\dok_m, P) > \delta\bigr) \leq C \exp
\biggl(\frac{-n\delta}{C^2 \sigma^2} \biggr) \qquad\mbox{whenever
$\delta
\geq\delta_0$}.
\]
Integrating both sides of the above inequality with respect to $\delta
\in[\delta_0, \infty)$, we get
\[
\E_K \bigl(\ell_f^2(\dok_m, P) - \delta_0 \bigr)_+ \leq
\frac{C\sigma^2}{n} \exp\biggl(\frac{-n\delta_0}{C^2\sigma
^2}\biggr)
\leq\frac{C\sigma^2}{n},
\]
where $x_+ := \max(x, 0)$. Because $\sigma^2/n \leq C \delta_0$, we get
the expectation bound $\E_K \ell_f^2(\dok_m, P) \leq C \delta_0$. The
elementary inequality $\ell_f^2(\dok_m, K) \leq2 \ell_f^2(\dok_m,
P) + 2
\ell_f^2(K, P)$ and the fact that $\delta_0 \geq\ell_f^2(K, P)$
yield the following risk bound:
%
\begin{equation}\label{rbmain}
\E_K \ell_f^2(\dok_m, K) \leq C \biggl( \ell_f^2(K, P) +
\frac{\sigma^2}{n}md \log(b_2m)\log\bigl(4+2\sqrt{n}\bigr) \biggr).
\end{equation}
Since $K \in\K^d(\Gamma)$ and $P \in\Ps_m(\Gamma)$ were arbitrary in
the above analysis, we have proved the following bound for
$\rok_f(n)$:
\[
\rok_f(n) \leq C \inf_{m \geq1} \biggl[\sup_{K \in\K^d(\Gamma)}
\inf_{P \in\Ps_m(\Gamma)} \ell_H^2(K, P) + \frac{\sigma^2}{n} md
\log
(b_2m)\log\bigl(4+2\sqrt{n}\bigr) \biggr],
\]
where we have also used $\ell_f^2(K, P) \leq\ell_H^2(K, P)$.\vspace*{2pt}

Bronshtein and
Ivanov~\cite{BronIvan} (see also Bronshtein~\cite{Bronshtein08}) proved that
there exist positive constants $C_1$ and $C_2$ depending only on the
dimension~$d$ such that
%
\begin{equation}\label{bronivaneq}\qquad
\sup_{K \in\K^d(\Gamma)} \inf_{P \in\Ps_m(\Gamma)} \ell_H^2(K, P)
\leq C_1 \Gamma^2 m^{-4/(d-1)} \qquad\mbox{whenever $m \geq C_2$}.
\end{equation}

From this result, we have
\[
\rok_f(n) \leq C \inf_{m \geq C_2} \biggl[\Gamma^2 m^{-4/(d-1)} +
\frac{\sigma^2}{n} md \log(b_2m)\log\bigl(4+2\sqrt{n}\bigr) \biggr].
\]
If we now choose $m := \sigma^{-2(d-1)/(d+3)}
\Gamma^{2(d-1)/(d+3)}n^{(d-1)/(d+3)}$, then $m \geq C_2$ provided $n
\geq C(\sigma/\Gamma)^2$ for a~large enough constant $C$ depending only
on~$d$ and the required inequality~(\ref{lsseq}) follows.
\end{pf*}
%
\begin{remark}
From the above proof, it can be seen that Theorem~\ref{lss} also
holds for the following estimator:
\[
\tilde{K}_m := \argmin_{L \in\Ps_m} \sum_{i=1}^n \bigl(Y_i -
h_L(u_i) \bigr)^2.
\]
The only difference between $\tilde{K}_m$ and $\dok_m$ is that the
argmin in the definition of $\tilde{K}_m$ is taken over all polytopes
in $\Ps_m$ while that in the definition of $\dok_m$ is only over those
polytopes in $\Ps_m$ that are contained in the ball of radius~$\Gamma$
centered at the origin. Theorem~\ref{lss} also holds for $\tilde{K}_m$
because we have not used this boundedness property of sets in
$\Ps_m(\Gamma)$ in our covering number calculations in Lemma~\ref{bpw}.
\end{remark}
%
\begin{remark}
The above proof also reveals that the value of $m$ for which
$\dok_m$ achieves the optimal rate up to logarithmic terms is of the
order $n^{(d-1)/(d+3)}$. Since this is much smaller than~$n$, the
estimator $\dok_m$ for this $m$ is quite different from the full
least squares estimator $\hat{K}_{\mathrm{ls}}$.
\end{remark}

\section{Random design setting}\label{rd}
In this section, we assume that $u_1,\ldots, u_n$ are independently
distributed according to a~fixed probability measure $\nu$ on the unit
sphere, $S^{d-1}$. The measurement errors $\xi_1,\ldots,
\xi_n$ are independent normal random variables with mean zero and
variance $\sigma^2$. We also assume that $\xi_1,\ldots, \xi_n$ are
independent of $u_1,\ldots, u_n$. We consider the loss function~$\ell_r^2$ [defined in~(\ref{rdlf})] and prove upper and lower bounds
for the corresponding minimax risk $R_r(n)$ [defined in~(\ref{rdmr})].

\subsection{Upper bound for $R_r(n)$}
The following result is the random-design analogue of
Theorem~\ref{ubfthm}. We show that $R_r(n)$ is bounded from above by
$n^{-4/(d+3)}$ up to multiplication by a~constant that is independent
of~$n$. No assumptions on $\nu$ are required and it is completely
arbitrary.
%
\begin{theorem}\label{ubrthm}
There exist positive constants $c$ and $C$ depending only on the
dimension~$d$ such that
%
\begin{equation}\label{ubrthmeq}
R_r(n) \leq\frac{c(\Gamma^2/\sigma^2)}{1 - e^{-\Gamma^2/(4
\sigma^2)}} \sigma^{8/(d+3)} \Gamma^{2(d-1)/(d+3)}
n^{-4/(d+3)},
\end{equation}
if $n \geq C(\sigma/\Gamma)^2$.
\end{theorem}
%
\begin{remark}\label{newlossrem}
Our proof below also shows that if one works with the smaller loss
function:
%
\begin{equation}\label{newloss}
\ell_{\mathrm{new}}^2(K, K') := - 16\sigma^2 \log\int\exp
\biggl(-\frac{(h_K(u) - h_{K'}(u))^2}{16\sigma^2} \biggr) \,d\nu(u)
\end{equation}
instead of $\ell_r^2$, then the factor
$(\Gamma^2/\sigma^2)(1-e^{-\Gamma^2/(4\sigma^2)})^{-1}$ in the minimax
risk bound~(\ref{ubrthmeq}) can be removed.
\end{remark}
\begin{pf*}{Proof of Theorem~\ref{ubrthm}}
As in the proof of Theorem~\ref{ubfthm}, we consi\-der the least
squares estimator $\hat{K}_F$ [defined in~(\ref{fls})] over a~finite
subset~$F$ of~$\K^d(\Gamma)$ for which inequality~(\ref{repineq}),
reproduced below, holds for every $K,\allowbreak L' \in F$ and $\alpha> 0$:
\[
\E_K G(\hat{K}_F) \leq\sum_{L \in F} G(L) \E_K \exp\biggl(\alpha
\sum_i \bigl(Y_i - h_{L'}(u_i)\bigr)^2 - \alpha\sum_i \bigl(Y_i - h_L(u_i)\bigr)^2
\biggr).
\]
Under the random-design setting, the expectation in the right-hand side
above equals
\[
\E_{\nu} \exp\bigl(-\alpha n \ell_f^2(K, L) + \alpha n \ell_f^2(K,
L') + 2\alpha^2 \sigma^2 n \ell_f^2(L, L') \bigr),
\]
where, as explained in Section~\ref{stanot}, the expectation
$\E_{\nu}$ is taken under the assumption that $u_1,\ldots, u_n$ are
independently distributed according to $\nu$ (note that
$\ell_f^2$ depends on $u_1,\ldots, u_n$). Using the inequalities
$\ell_f^2(L, L') \leq2 \ell_f^2(K, L) + 2 \ell_f^2(K, L')$ and
$\ell_f^2(K, L') \leq\ell_H^2(K, L')$, we obtain the
following upper bound for $\E_K G(\hat{K}_F)$:
\[
\sum_{L \in F} \exp\bigl((\alpha+ 4 \alpha^2
\sigma^2) n \ell_H^2(K, L') \bigr) G(L) \E_{\nu} \exp
\bigl((-\alpha+ 4 \alpha^2 \sigma^2) n \ell_f^2(K, L) \bigr)
\]
for every $L' \in F$. It may be helpful to note that $\ell_H^2(K, L')$
does not depend on $u_1,\ldots, u_n$ and is nonrandom. We apply this
inequality to the following choices of $G$ and $\alpha$:
\[
G(L) = \biggl(\E_{\nu} \exp\biggl(- \frac{n \ell_f^2(K, L)}{16
\sigma^2} \biggr) \biggr)^{-1}  \quad\mbox{and}\quad \alpha=
\frac{1}{8 \sigma^2}.
\]
A~straightforward calculation reveals that this function $G(L)$ has
the following alternative expression:
\[
G(L) = \exp\biggl(\frac{n}{16\sigma^2} \ell^2_{\mathrm{new}}(K, L) \biggr),
\]
where $\ell^2_{\mathrm{new}}$ is defined as in~(\ref{newloss}). Specializing
the upper bound for $\E_KG(\hat{K}_F)$ to these choices of $G$ and
$\alpha$, we deduce
\[
\E_K \exp\biggl(\frac{n \ell_{\mathrm{new}}^2(K, \hat{K}_F)}{16 \sigma^2}
\biggr)
\leq|F| \exp\biggl(\frac{3n}{16 \sigma^2} \min_{L' \in F}
\ell_H^2(K, L') \biggr).
\]
This is the same inequality as~(\ref{expb}) with the loss function
$\ell^2_f$ replaced by $\ell^2_{\mathrm{new}}$. Thus, following the same steps
as in the proof of Theorem~\ref{ubfthm}, we deduce the existence of
positive constants $c$ and $C$ depending only on~$d$ such that
%
\begin{equation}\label{newmm}
R_{\mathrm{new}}(n) \leq c \sigma^{8/(d+3)} \Gamma^{2(d-1)/(d+3)}n^{-4/(d+3)}
\qquad\mbox{if $n \geq C(\sigma/\Gamma)^2$},
\end{equation}
where
\[
R_{\mathrm{new}}(n) := \inf_{\hat{K}} \sup_{K \in\K^d(\Gamma)} \E_K
\ell_{\mathrm{new}}^2(K, \hat{K}).
\]
This proves the claim made in Remark~\ref{newlossrem}. We now give a~simple inequality relating the loss functions $\ell_{\mathrm{new}}^2(K, K')$
and $\ell_{r}^2(K, K')$ which enables us to convert this bound for
$R_{\mathrm{new}}(n)$ into the required inequality~(\ref{ubrthmeq}) for
$R_r(n)$. For every $K \in\K^d(\Gamma)$ and $u \in S^{d-1}$, we have
$|h_K(u)| \leq\Gamma$. Therefore,
\[
\frac{(h_K(u) - h_{K'}(u))^2}{16 \sigma^2} \leq\frac{\Gamma^2}{4
\sigma^2} \qquad\mbox{for all $K, K' \in\K^d(\Gamma)$}.
\]
Since the convex function $x \mapsto e^{-x}$ lies below the chord
joining the points $(0, 1)$ and $(\Gamma^2/(4 \sigma^2), \exp
(-\Gamma^2/(4 \sigma^2)))$, we have
\[
e^{-x} \leq1 + \frac{4\sigma^2}{\Gamma^2} \bigl(e^{-\Gamma^2/(4
\sigma^2)} - 1 \bigr)x \qquad\mbox{for $0 \leq x \leq
\Gamma^2/(4\sigma^2)$}.
\]
Using this with $x = (h_K(u) - h_{K'}(u))^2/(16\sigma^2)$ and
integrating both sides of the resulting expression with respect to
$\nu$, we get
\[
\int\exp\biggl(- \frac{(h_K(u) - h_{K'}(u))^2}{16 \sigma^2} \biggr)
\,d\nu(u) \leq1 + \frac{e^{-\Gamma^2/(4\sigma^2)} - 1}{4\Gamma^2}
\ell_r^2(K, K').
\]
Taking logarithms on both sides, we obtain
\[
\log\int\exp\biggl(- \frac{(h_K(u) - h_{K'}(u))^2}{16 \sigma^2}
\biggr)
\,d\nu(u) \leq\frac{e^{-\Gamma^2/(4\sigma^2)} - 1}{4\Gamma^2}
\ell_r^2(K, K'),
\]
where, on the right-hand side, we have used $\log(1 + y) \leq y$. The
above inequality can be rewritten as
\[
\ell_r^2(K, K') \leq\frac{\Gamma^2/(4\sigma^2)}{1 -
e^{-\Gamma^2/(4\sigma^2)}} \ell^2_{\mathrm{new}}(K, K').
\]
The proof is complete because the required bound~(\ref{ubrthmeq})
follows by combining the above inequality with~(\ref{newmm}).
\end{pf*}

\subsection{Lower bound for $R_r(n)$}
The following theorem is the random-design analogue of
Theorem~\ref{assf}. We assume that $\nu= \unisp$ is the uniform
probability measure on $S^{d-1}$ and prove that $R_r(n)$ is bounded
from below by a~constant multiple of $n^{-4/(d+3)}$. Note that the
lower bound of $n^{-4/(d+3)}$ cannot be proved for $R_r(n)$ for
arbitrary $\nu$. For example, when $\nu$ is concentrated at a~single
point, $R_r(n)$~is of order $1/n$.
%
\begin{theorem}\label{lbrthm}
Consider the random-design setting where $\nu$ equals the uniform
probability measure $\unisp$ on $S^{d-1}$. Then there exist positive
constants $c$ and $C$ depending only on~$d$ such that
%
\begin{equation}
\label{lbrthmeq}
R_r(n) \geq c \sigma^{8/(d+3)} \Gamma^{2(d-1)/(d+3)} n^{-4/(d+3)}
\qquad\mbox{whenever}\ n \geq C (\sigma/\Gamma)^2.\hspace*{-30pt}
\end{equation}
\end{theorem}

It is possible to prove this theorem by an appropriate modification of
the proof of Theorem~\ref{assf}. We, however, give a~different proof
using a~global metric entropy minimax lower bound
from Guntuboyina~\cite{GuntuFdiv}. This proof has an interesting implication
that is described in Remark~\ref{arbdes}. A~version of this proof
appeared in Guntuboyina~\cite{GuntuFdiv}, Section V, although the
result there has
a~different assumption on $u_1,\ldots, u_n$ and is also slightly less
precise.
\begin{pf*}{Proof of Theorem~\ref{lbrthm}}
Let $\mathfrak{P} := \{ \P_K \dvtx K \in\K^d(\Gamma)\}$. We use the
following minimax lower bound from Guntuboyina~\cite{GuntuFdiv},
inequality (22):
%
\begin{equation}\label{sless}
R_r(n) \geq\frac{\eta^2}{4} \Biggl(1
- \frac{1}{N(\K^d(\Gamma), \eta; \ell_r)} -
\sqrt{\frac{(1+\epsilon^2)M(\mathfrak{P}, \epsilon;
\chi)}{N(\K^d(\Gamma), \eta; \ell_r)}}\Biggr)
\end{equation}
for all $\eta> 0$ and $\epsilon> 0$. The notation was set in
Section~\ref{stanot}: $N(\K^d(\Gamma), \eta; \ell_r)$ denotes the
maximal $\eta$-packing number of $\K^d(\Gamma)$ under the
$\ell_r$-metric (the square root of the loss function $\ell_r^2$) and
$M(\mathfrak{P}, \epsilon; \chi)$ denotes the $\epsilon$-covering
number of $\mathfrak{P}$ when distances are measured by the
square root of the chi-squared divergence, that is, it is the smallest
integer $M$ for which there exist probability measures $\Q_1,\ldots,
\Q_M$ satisfying $\min_{1 \leq i \leq M} \chi^2(\P\vvert\Q_i) \leq
\epsilon^2$ for every $\P\in\mathfrak{P}$.

The application of~(\ref{sless}) requires a~lower bound on
$N(\K^d(\Gamma), \eta; \ell_r)$ and an upper bound on $M(\mathfrak{P},
\epsilon; \chi)$. Guntuboyina~\cite{GuntuFdiv}, Theorem VII.1,
building on a~result of Bronshtein~\cite{Bronshtein76}, showed the
existence of positive
constants $\eta_0$ and~$c'$ depending only on~$d$ such that
\[
\log N(\K^d(\Gamma), \eta; \ell_r) \geq c'
\biggl(\frac{\Gamma}{\eta} \biggr)^{(d-1)/2} \qquad\mbox{whenever
$\eta\leq\Gamma\eta_0$}.
\]
The above bound uses crucially the fact that $\nu$ equals $\unisp$. It
is not true for arbitrary probability measures on $S^{d-1}$.

For $M(\mathfrak{P}, \epsilon; \chi)$, we note that the chi-squared
divergence $\chi^2(\P_K\vvert\P_{K'})$ satisfies
\begin{eqnarray*}
1+\chi^2(\P_K\vvert\P_{K'}) &=& \biggl(\int\exp\biggl(\frac{(h_K(u) -
h_{K'}(u))^2}{\sigma^2} \biggr) \,d\nu(u) \biggr)^n \\
&\leq& \exp\biggl(\frac{n\ell_H^2(K, K')}{\sigma^2} \biggr).
\end{eqnarray*}
As a~result,
\[
\chi^2(\P_K\vvert\P_{K'}) \leq\epsilon^2 \qquad\mbox{whenever $\ell_H(K,
K') \leq
\epsilon' := \sigma\sqrt{\log(1+\epsilon^2)}/\sqrt{n}$}
\]
and $M(\mathfrak{P}, \epsilon; \chi) \leq M(\K^d(\Gamma), \epsilon';
\ell_H)$. Upper bound for the covering number $M(\K^d(\Gamma),
\epsilon';
\ell_H)$ has been proved by Bronshtein~\cite{Bronshtein76}, Theorem 3 and Remark~1. We already used this result
[inequality~(\ref{bron1})] in the proof of Theorem~\ref{ubfthm}. We
use it again to obtain
\[
\log M(\mathfrak{P}, \epsilon; \chi) \leq c \biggl(\frac{\Gamma
\sqrt{n}}{\sigma\sqrt{\log(1+\epsilon^2)}} \biggr)^{(d-1)/2}
\qquad\mbox{if $\log(1+\epsilon^2) \leq n \Gamma^2
\epsilon_0^2/\sigma^2$}
\]
for positive constants $c$ and $\epsilon_0$.
Let us now define
\[
\eta(n) := c_1 \sigma^{4/(d+3)} \Gamma^{(d-1)/(d+3)}n^{-2/(d+3)}
\quad\mbox{and}\quad \alpha(n) := \biggl(\frac{\Gamma\sqrt{n}}{\sigma}
\biggr)^{(d-1)/(d+3)},
\]
where $c_1$ is a~positive constant that depends on~$d$ alone and
will be specified shortly. Also let $\epsilon^2(n) :=
\exp(\alpha^2(n))-1$. We then have
\[
\log N(\K^d(\Gamma), \eta(n); \ell_r) \geq c' c_1^{-(d-1)/2}
\alpha^2(n) \quad\mbox{and}\quad\log M(\mathfrak{P}, \epsilon(n); \chi)
\leq c\alpha^2(n)
\]
provided
%
\begin{equation}\label{reqs}
\eta(n) \leq\Gamma\eta_0 \quad\mbox{and}\quad\alpha^2(n) \leq n
\Gamma^2 \epsilon_0^2/\sigma^2.
\end{equation}
Inequality~(\ref{sless}) with $\eta= \eta(n)$ and $\epsilon=
\epsilon(n)$ gives the following lower bound for $R_r(n)$:
\[
\frac{\eta^2(n)}{4} \biggl[1 - \exp\bigl(-\alpha^2(n) c'c_1^{-(d-1)/2}
\bigr) - \exp\biggl(\frac{\alpha^2(n)}{2}\bigl(1+c - c'c_1^{-(d-1)/2}\bigr)
\biggr) \biggr].
\]
If we choose $c_1$ so that $c'c_1^{-(d-1)/2} = 2(1+c)$, then
\[
R_r(n) \geq\frac{\eta^2(n)}{4} \biggl(1 - 2 \exp
\biggl(-\frac{1+c}{2}\alpha^2(n) \biggr) \biggr).
\]
If the condition $(1+c)\alpha^2(n) \geq2\log4$ holds, then the
above inequality implies $R_r(n) \geq\eta^2(n)/8$ which
yields~(\ref{lbrthmeq}). This condition as well as~(\ref{reqs}) hold
provided $n \geq C(\sigma/\Gamma)^2$ for a~large enough $C$. The proof
is complete.
\end{pf*}
%
\begin{remark}\label{arbdes}
In the above proof, the random-design assumption on the unit vectors
$u_1,\ldots, u_n$ was used only in
\[
\chi^2(\P_K\vvert\P_{K'}) \leq\exp\biggl(\frac{n\ell_H^2(K,
K')}{\sigma^2} \biggr) - 1.
\]
This inequality is easily seen to be true for every joint distribution
of $(u_1,\ldots,\allowbreak u_n)$ as long as they are independent of the errors
$\xi_1,\ldots, \xi_n$. Consequently, $n^{-4/(d+3)}$, up to
multiplicative constants, is a~lower bound for the minimax risk
(observe that the integral below is with respect to the uniform
probability measure~$\nu_{\mathrm{unif}}$):
\[
\inf_{\hat{K}} \sup_{K \in\K^d(\Gamma)} \E_K \int_{S^{d-1}}
\bigl(
h_K(u) - h_{\hat{K}}(u) \bigr)^2 \unisp(du)
\]
for every arbitrary choice of the design unit vectors (deterministic
or random) as long as they are independent of $\xi_1,\ldots, \xi_n$
(this independence assumption is only relevant in a~random-design
setting).
\end{remark}

\subsection{Least squares on polytopes}
We prove the random-design analogue of Theorem~\ref{lss}. Let
\[
\rok_r(n) := \inf_{m \geq1} \sup_{K \in\K^d(\Gamma)} \E_K
\ell_r^2(K, \dok_m)
\]
be the best achievable risk by $\dok_m$ [defined in~(\ref{dokyun})] in
the random-design setting under the loss function $\ell_r^2$. The
following theorem shows that $\rok_r(n)$ is bounded from above by
$n^{-4/(d+3)}$ up to a~multiplicative factor that is logarithmic in
$n$. No assumptions on $\nu$ are necessary.
%
\begin{theorem}\label{rlss}
There exist positive constants $c$ and $C$ that depend only on the
dimension~$d$ such that
%
\begin{equation}
\label{rlsseq}
\rok_r(n) \leq c \max(\sigma^2, \Gamma^2) n^{-4/(d+3)} (
\log(cn) )^2 \qquad\mbox{if $n \geq C$}.
\end{equation}
\end{theorem}

The proof strategy is to use the fixed design bound~(\ref{lsseq})
along with Lem\-ma~\ref{gyopol} (stated and proved in Appendix~\ref{appC})
which relates the risks under the two loss functions $\ell_f^2$ and
$\ell_r^2$.
\begin{pf*}{Proof of Theorem~\ref{rlss}}
We start with the inequality
\[
\ell_r^2(\dok_m, K) \leq2 \bigl(\ell_r(\dok_m, K) - 2
\ell_f(\dok_m, K) \bigr)_+^2 + 8 \ell_f^2(\dok_m, K),
\]
which implies that
\[
\E_K \ell_r^2(\dok_m, K) \leq2 \E_{\nu} \sup_{L \in
\Ps_m(\Gamma)} \bigl(\ell_r(L, K) - 2
\ell_f(L, K) \bigr)_+^2 + 8 \E_K \ell_f^2(\dok_m, K).
\]
The first expectation in the right-hand side is bounded using
Lemma~\ref{gyopol} where it is shown that
%
\begin{equation}\label{pbook}
\E_{\nu} \sup_{L \in\Ps_m(\Gamma)} \bigl(\ell_r(L, K) - 2
\ell_f(L, K) \bigr)_+^2 \leq c \frac{\Gamma^2}{n} md\log(cn)
\end{equation}
for a~universal positive constant $c$. For the second expectation, we
use ideas from the proof of Theorem~\ref{lss}. Indeed, the same
argument which led to the inequality~(\ref{rbmain}) gives, for every
$P \in\Ps_m$,
\[
\E_K (\ell_f^2(\dok_m, K) \vert u_1,\ldots, u_n ) \leq c
\biggl( \ell_f^2(K, P) + \frac{\sigma^2}{n}md \log
(b_2m)\log(cn) \biggr)
\]
for every $u_1,\ldots, u_n$. Taking expectation with respect to $u_1,\ldots, u_n$ independently distributed according to $\nu$, we get
%
\begin{equation}\label{rbagain}
\E_K \ell_f^2(\dok_m, K) \leq c \biggl( \ell_r^2(K, P) +
\frac{\sigma^2}{n}md \log(b_2m)\log(cn) \biggr).
\end{equation}
Putting~(\ref{pbook}) and~(\ref{rbagain}) together, we obtain
\begin{eqnarray*}
\E_K \ell_r^2(\dok_m, K ) &\leq& c \biggl(\ell_r^2(K,
P) +
\frac{\sigma^2}{n}md \log(b_2m)\log(cn) +
\frac{\Gamma^2}{n} md \log(cn) \biggr) \\
&\leq& c \biggl(\ell_H^2(K, P) +
\frac{\max(\sigma^2, \Gamma^2)}{n}md \log(b_2m)\log(cn) \biggr).
\end{eqnarray*}
Because $K \in\K^d(\Gamma)$ and $P \in\Ps_m(\Gamma)$ are arbitrary,
we have shown
\[
\rok_r(n) \leq c \inf_{m \geq1} \biggl[\sup_{K \in\K^d(\Gamma)}
\inf_{P \in\Ps_m(\Gamma)} \ell_H^2(K, P) + \frac{\max(\sigma^2,
\Gamma^2)}{n}md \log(b_2m)\log(cn) \biggr].
\]
Just as in the proof of Theorem~\ref{lss}, we use the
result~(\ref{bronivaneq}) due to Bronshtein and
Ivanov~\cite{BronIvan} to get
\begin{eqnarray*}
\rok_r(n) &\leq& C \inf_{m \geq C_2} \biggl[\Gamma^2 m^{-4/(d-1)} +
\frac{\max(\sigma^2, \Gamma^2)}{n} md \log(b_2m)\log(cn) \biggr]
\\
&\leq& C \max(\sigma^2, \Gamma^2) \inf_{m \geq C_2} \biggl[ m^{-4/(d-1)}
+ \frac{md}{n} \log(b_2m)\log(cn) \biggr].
\end{eqnarray*}
If we now choose $m := n^{(d-1)/(d+3)}$, then $m \geq C_2$ provided $n
\geq C$ for a~large enough constant $C$ depending only
on~$d$ and the required inequality~(\ref{rlsseq})
follows.\hspace*{0pt}\qed
\noqed\end{pf*}
\begin{appendix}\label{app}
%
\section{}\label{appA}
In this section, we shall prove the following result which was used
in the proof of Theorem~\ref{assf}. We assume that $u_1,\ldots, u_n$
form a~maximal $\epsilon$-packing subset of $S^{d-1}$ for some
$\epsilon\in(0, 1]$.
%
\begin{lemma}\label{appass}
For a~fixed $0 < \eta\leq1/8$ and a~unit vector $v$, consider the
following two subsets of the ball $B_d(0, \Gamma)$:
\[
D(0) := B_d(0, \Gamma) \cap\{x \dvtx x
\cdot v \leq1 - \eta\}  \quad\mbox{and}\quad  D(1) :=
B_d(0, \Gamma).
\]
Then there exists constants $c$ and $C$ such that the following inequality
holds whenever $\eta\geq C \epsilon^2$:
%
\begin{equation}\label{l2capsupp}
c \eta^{(d+3)/2} \leq\ell_f^2(D(0), D(1)) \leq
C \eta^{(d+3)/2}.
\end{equation}
\end{lemma}

We need some elementary results on spherical caps for the proof of this
lemma. For a~unit vector $u$ and a~real number $0 < \delta\leq
1$, consider the spherical cap $\scap(u;\delta) := S^{d-1} \cap B_d(u,
\delta)$. It can be checked that this spherical cap consists of
precisely those unit vectors which form an angle of at most $\alpha$
with $u$, where $\alpha$ is related to $\delta$ through
\[
\cos\alpha= 1 - \frac{\delta^2}{2} \quad\mbox{and}\quad \sin\alpha=
\frac{\delta\sqrt{4 - \delta^2}}{2}.\vadjust{\goodbreak}
\]
A~standard result is that $\unisp(\scap(x; \delta))$ equals $C
\int_0^{\alpha} \sin^{d-2}t \,dt$ (recall that~$\unisp$ denotes the uniform
probability measure on the unit sphere). This integral can be bounded
from above and below in the following simple way. For a~lower bound,
we write
\[
\int_0^{\alpha} \sin^{d-2}t \,dt \geq\int_0^{\alpha} \sin^{d-2}t
\cos t \,dt \geq\frac{\sin^{d-1} \alpha}{d-1},
\]
and for an upper bound, we note
\[
\int_0^{\alpha} \sin^{d-2}t \,dt \leq\int_0^{\alpha}
\frac{\cos t}{\cos\alpha} \sin^{d-2}t \,dt \leq\frac{\sin^{d-1}
\alpha}{(d-1)\cos\alpha}.
\]
We thus have $c \sin^{d-1} \alpha\leq\nu(\scap(u; \delta)) \leq C
\sin^{d-1} \alpha/\cos\alpha$. Writing $\cos\alpha$ and $\sin
\alpha$ in terms of $\delta$ and using the assumption that $0 <
\delta
\leq1$, we obtain that
%
\begin{equation}\label{mucap1}
c \delta^{d-1} \leq\unisp(\scap(u; \delta)) \leq C \delta^{d-1}
\qquad\mbox{if
$\delta\in(0, 1]$}.
\end{equation}
This inequality can be combined with a~simple
volumetric argument to show that
%
\begin{equation}\label{pent}
c \delta^{1-d} \leq N(S^{d-1}, \delta; \euc_d) \leq C \delta^{1-d}
\qquad\mbox{if $\delta\in(0, 1]$.}
\end{equation}
In particular, since $u_1,\ldots, u_n$ form a~maximal
$\epsilon$-packing subset of $S^{d-1}$, we
have~(\ref{maxpack}).

The following lemma is used in the proof of Lemma~\ref{appass}.
%
\begin{lemma}\label{vep}
Fix positive $\epsilon, \delta$ such that $\delta\leq4/5$ and
$\epsilon\leq\delta/2$. Let $u_1, \ldots, u_n$ be a~maximal
$\epsilon$-packing subset of the unit sphere and $v$ be an arbitrary
unit vector. Let $V(\epsilon, \delta)$ denote the number
of points $u_1,\ldots, u_n$ that are contained in $\scap(v,
\delta)$. Then
%
\begin{equation}\label{vb}
c \biggl(\frac{\delta}{\epsilon}\biggr)^{d-1} \leq V(\epsilon,
\delta)
\leq C \biggl(\frac{\delta}{\epsilon}\biggr)^{d-1}.
\end{equation}
\end{lemma}
\begin{pf}
For the lower bound on $V(\epsilon, \delta)$, we observe that
%
\begin{equation}\label{bageshri}
\scap(v, \delta/2) \subseteq\bigcup_{i\dvtx u_i \in\scap(v, \delta)}
\scap(u_i, \epsilon).
\end{equation}
Indeed, because $\epsilon\leq\delta/2$, for every $w \in\scap(v,
\delta/2)$, we can find $u_i$ such that $\euc_d(u_i, w) < \epsilon
\leq
\delta/2$ because $u_1,\ldots, u_n$ form a~maximal $\epsilon$-packing
subset\break of $S^{d-1}$. Thus $u_i \in\scap(v, \delta)$ by triangle
inequality which proves~(\ref{bageshri}).

It follows from~(\ref{bageshri}) that $V(\epsilon, \delta) \geq
\unisp(\scap(v, \delta/2))/\unisp(\scap(u_1, \epsilon))$ from
which the
lower bound in~(\ref{vb}) follows by use of~(\ref{mucap1}).

For the upper bound on $V(\epsilon, \delta)$, we use the inequality
\[
\scap(v, \delta+ \epsilon/2) \supseteq\bigcup_{i\dvtx u_i \in\scap(v,
\delta)} \scap(u_i, \epsilon/2)
\]
and then, noting that the spherical caps on the right-hand side of the
above inequality are disjoint as $u_1,\ldots, u_n$ form a~$\epsilon$-packing subset, we obtain $V(\epsilon, \delta) \leq
\unisp(\scap(v, \delta+
\epsilon/2))/\unisp(\scap(u_1,\epsilon/2))$. The upper bound
in~(\ref{vb}) again follows from~(\ref{mucap1}).
\end{pf}

We are now ready to prove Lemma~\ref{appass}.
\begin{pf*}{Proof of Lemma~\ref{appass}}
It can be checked that the support functions of $D(0)$ and $D(1)$
differ only for unit vectors in the spherical cap
$\scap(v, \sqrt{2 \eta})$. This spherical cap consists of all
unit vectors which form an angle of at most $\alpha$ with
$v$ where $\cos\alpha= 1 - \eta$. In fact, if $\theta$ denotes the
angle between an arbitrary unit vector $u$ and~$v$, it can be
verified by elementary trigonometry that
%
\begin{equation}\label{exsupp}
h_{D(0)}(u) - h_{D(1)}(u) =
\cases{
\Gamma\bigl(1 - \cos(\alpha- \theta
)\bigr), &\quad if $0 \leq\theta\leq\alpha$,\cr
0, &\quad otherwise.}
\end{equation}
As a~result, it follows that $\ell_f^2(D(0), D(1))$ $\leq$
$\Gamma^2 \eta^2 V(\epsilon, \sqrt{2\eta})/n$ where $V$ is
as defined in Lemma~\ref{vep}. Thus~(\ref{vb}) gives
$\ell_f^2(D(0), D(1))$ $\leq C \eta^{(d+3)/2}\epsilon^{1-d}/n$. The
conditions in Lemma~\ref{vep} are satisfied if $0 < \eta\leq1/8$ and
$\eta\geq C \epsilon^2$ for a~large enough $C$. Moreover,
by~(\ref{maxpack}), we have $c \leq n\epsilon^{d-1} \leq C$ which
implies that $\ell_f^2(D(0), D(1)) \leq C \eta^{(d+3)/2}$.

For a~lower bound, fix $0<b\leq1$ and let $0 \leq\beta\leq\alpha$
denote the angle for which $1 - \cos(\alpha- \beta) = b\eta$. It follows
from~(\ref{exsupp}) that the difference in the support functions of
$D(0)$ and $D(1)$ is at least $b \Gamma\eta$ for all unit vectors
in the spherical cap consisting of all unit vectors forming an angle
of at most $\beta$ with~$v$. This spherical cap is
$\scap(v, t)$ where~$t$ is given by $t^2 := 2(1-\cos
\beta)$. Therefore $\ell_f^2(D(0), D(1)) \geq b^2 \Gamma^2
\eta^2 V(\epsilon, t)/n$. The\vspace*{1pt} inequality $t^2 \leq2(1-\cos
\alpha) \leq2 \eta$ is easily checked. Also, $t \geq\sin\beta$
and $\sin\beta$ can be bounded from below in the following way:
\[
1-b\eta= \cos(\alpha- \beta) \leq\cos\alpha+ \sin\alpha\sin
\beta\leq1-\eta+ \sqrt{2\eta} \sin\beta.
\]
Thus $t \geq\sin\beta\geq(1-b) \sqrt{\eta/2}$ and
from~(\ref{vb}), it follows that
\[
\ell_f^2(D(0), D(1)) \geq c b^2 \frac{\Gamma^2
\eta^2}{n}\biggl(\frac{t}{\epsilon}\biggr)^{d-1} \geq c b^2
(1-b)^{d-1} \frac{\Gamma^2\eta^{(d+3)/2}}{n\epsilon^{d-1}}
\]
for all $0 < b \leq1$. Note that we have used $\eta
\geq C \epsilon^2$ here to satisfy the conditions in
Lemma~\ref{vep}. We now use~(\ref{maxpack}) and choose $b = 1/2$ to
get $\ell_f^2(D(0),\allowbreak D(1)) \geq c \eta^{(d+3)/2}$ provided $\eta
\geq C \epsilon^2$. The proof is complete.
\end{pf*}

\section{}\label{appB}
In this section, we prove the following lemma which was used in the
proof of Theorem~\ref{lss}.\vadjust{\goodbreak}
%
\begin{lemma}\label{bpw}
Fix $m \geq1$ and $\omega, \epsilon> 0$. The following bound holds
for every $P \in\Ps_m$:
\[
M(S_m(P, \omega), \epsilon; \ell_f) \leq\biggl(4 +
\frac{2\sqrt{n}\omega}{\epsilon} \biggr)^{b_1 md\log(b_2 m)},
\]
where $b_1$ and $b_2$ are universal positive constants.
\end{lemma}

For the proof of this lemma, we use available techniques for bounding
covering numbers using combinatorial notions of
dimension. Specifically, we use the notion of
pseudodimension, introduced by Pollard~\cite{Pollard90Iowa}, Chapter
4, as a~generalization of the Vapnik--\u{C}ervonenkis dimension to classes of
real-valued functions. The pseudodimension of a~subset $A$ of $\R^n$
is defined as the maximum cardinality of a~subset
$\sigma\subseteq\{1,\ldots, n\}$ for which there exists $(h_1,\ldots, h_n) \in\R^n$ such that: for every $\sigma' \subseteq\sigma$,
one can find $(a_1,\ldots, a_n) \in A$ with $a_i < h_i$ for $i \in
\sigma'$ and $a_i > h_i$ for $i \in\sigma\setminus\sigma'$. The
following theorem is a~special case of results in Pollard~\cite{Pollard90Iowa},
Chapter 4, and gives an upper bound for the covering
number (with respect to the Euclidean metric) of a~subset of $\R^n$ in
terms of its pseudodimension. Stronger results of this kind have been
proved by Mendelson and
Vershynin~\cite{MendelsonVershynin03} and the following theorem
is also a~special case of Mendelson and
Vershynin~\cite{MendelsonVershynin03}, Theorem~1.
%
\begin{theorem}\label{mvpd}
Let $A$ be a~subset of $\R^n$ with ${\max_i }|a_i| \leq B$ for all $a~=
(a_1,\ldots, a_n) \in A$. If the pseudodimension of $A$ is at most
$V$, then, for every $t > 0$, we have
\[
M(A, t; \euc_n) \leq\biggl(4 + \frac{2B\sqrt{n}}{t} \biggr)^{b V},
\]
where $b$ is a~universal positive constant.
\end{theorem}
\begin{pf*}{Proof of Lemma~\ref{bpw}}
Fix a~polytope $P \in\Ps_m$ and let $x_0 \in\R^n$ denote the point
$(h_P(u_1),\ldots, h_P(u_n))$. Also, for $m \geq1$, let
\[
H_m := \{x \in\R^n \dvtx x = (h_L(u_1),\ldots, h_L(u_n)) \mbox{
for some } L \in\Ps_m \}.
\]
Clearly
\[
M(S_m(P, \omega),\epsilon; \ell_f) = M\bigl(B_n\bigl(0, \sqrt{n}\omega\bigr)
\cap H_m - x_0, \sqrt{n}\epsilon; \euc_n\bigr),
\]
where $H_m - x_0 := \{ x - x_0 \dvtx x \in H_m\}$.

We now show that the pseudodimension of $H_m$, which clearly equals
the pseudodimension of $H_m - x_0$, is less than or equal to
$\varsigma md\log(\varsigma m)$ where $\varsigma= 2/\log2$. An
application of Theorem~\ref{mvpd} would then complete the proof. Note
that the quantity $B$ in the statement of Theorem~\ref{mvpd} can be
taken to be $\sqrt{n}\omega$ because $|a_i| \leq\sqrt{n}\omega$ for
every $(a_1,\ldots, a_n) \in B_n(0, \sqrt{n}\omega)$.

Every $L$ in $\Ps_m$ is of the form $\operatorname{conv}(S)$ for some $S
\subset\R^d$ with cardinality at most $m$ and thus $h_L(u) = \max_{x
\in S} (x \cdot u)$. Therefore $H_1$ is a~linear\vadjust{\goodbreak}
subspace of $\R^n$ with dimension at most~$d$ which implies
(see Pollard~\cite{Pollard90Iowa}, page 15) that the pseudodimension
of $H_1$
is at most~$d$. The pseudodimension of $H_m$, which consists of
coordinatewise maxima of at most $m$ points in $H_1$, can then be
bounded from above following the argument in Pollard~\cite{Pollard90Iowa},
proof of Lemma (5.1). Indeed, that argument shows that the
pseudodimension of~$H_m$ is bounded from above by the smallest
positive integer $k$ for which
%
\begin{equation}\label{reqvcd}
\pmatrix{k \cr 0} +\cdots+ \pmatrix{k \cr d} < 2^{k/m}.
\end{equation}
If $\mathfrak{B}$ denotes a~binomial random variable with parameters
$k$ and $1/2$, then the left-hand side of the above inequality equals
$2^k \P\{\mathfrak{B} \geq k-d \}$ and can be bounded,
again following Pollard~\cite{Pollard90Iowa}, proof of Lemma (5.1), as shown
below. For every $\alpha> 0$,
%
\begin{equation}\label{bino}\qquad
2^k \P\{\mathfrak{B} \geq k-d \} \leq2^k \P
\{\alpha^{\mathfrak{B}} \geq\alpha^{k-d} \} \geq2^k \alpha^{d-k}
\E\alpha^{\mathfrak{B}} = (1+\alpha)^k \alpha^{d-k}.
\end{equation}
If we now choose $\alpha= (2^{1/(2m)} - 1)^{-1}$, then
$(1+\alpha)/\alpha= 2^{1/(2m)}$ and also, applying the inequality $x
- 1 > \log x$ to $x = 2^{1/(2m)}$, we get that $\alpha< 2m/(\log
2)$. Therefore, from~(\ref{bino}), we have
\[
2^k \P\{\mathfrak{B} \geq k-d \} < \biggl(\frac{2m}{\log
2}\biggr)^d 2^{k/(2m)}.
\]
The following inequality therefore ensures that $k$
satisfies~(\ref{reqvcd}):
\[
\biggl(\frac{2m}{\log2}\biggr)^d 2^{k/(2m)} \leq2^{k/m} \quad\mbox{or,
equivalently}\quad k \geq\varsigma md\log(\varsigma m),
\]
where $\varsigma= 2/(\log2)$. It follows, therefore, that the
pseudodimension of $H_m$ is at most $\varsigma md\log(\varsigma
m)$. The proof of Lemma~\ref{bpw} is complete.
\end{pf*}

\section{}\label{appC}
In this section, we prove the following result which was used in the
proof of Theorem~\ref{rlss}. It relates the risks under the loss
functions $\ell_f^2$ and $\ell_r^2$. Its proof is based on standard
empirical process arguments due to Pollard~\cite{Pollard84book}. We
have also
borrowed ideas from
Gy{\"o}rfi, Kohler, Krzy{\.z}ak and Walk
\cite{Gyorfi}, Chapter 11.
%
\begin{lemma}\label{gyopol}
There exists a~universal constant $c$ such that the following
inequality holds for every $m \geq1$, $\Gamma> 0$ and $K \in
\Ps_m(\Gamma)$:
%
\begin{equation}\label{gyopoleq}
\E_{\nu} \sup_{L \in\Ps_m(\Gamma)} \bigl(\ell_r(L, K) - 2
\ell_f(L, K) \bigr)_+^2 \leq c \frac{\Gamma^2}{n} md\log(cn).
\end{equation}
\end{lemma}
\begin{pf}
For $x > 0$, let
\[
\varrho_x := \P_{\nu} \Bigl\{\sup_{L \in\Ps_m(\Gamma)}
\bigl(\ell_r(L, K) - 2 \ell_f(L, K) \bigr)_+^2 > x \Bigr\}.\vadjust{\goodbreak}
\]
Letting $\F:= \{h_L - h_K \dvtx L \in\Ps_m(\Gamma)\}$, we have the
trivial inequality
\[
\varrho_x \leq\P_{\nu} \Biggl\{\sup_{f \in\F} \Biggl[\biggl(\int f^2(u)
\,d\nu(u)\biggr)^{1/2} - 2 \Biggl(\frac{1}{n} \sum_{i=1}^n f^2(u_i)
\Biggr)^{1/2}\Biggr] > \sqrt{x} \Biggr\}.
\]
Each function in $\F$ is bounded in absolute value by $2
\Gamma$. We now use Gy{\"o}rfi et al.~\cite{Gyorfi}, Theorem 11.2,
which is a~slight
variation of Pollard~\cite{Pollard84book}, Lemma~33 and Problem 24,
to obtain
%
\begin{equation}\label{akka}
\varrho_x \leq4 \E_{\nu} \min\biggl(1, M\bigl(\Ps_m(\Gamma), \sqrt{x}/24;
\ell_f\bigr) \exp\biggl(-\frac{nx}{2304 \Gamma^2} \biggr) \biggr).
\end{equation}
Lemma~\ref{bpw} with $\omega= \Gamma$, $\epsilon= \sqrt{x}/24$ and
$P = \{0\}$ provides an upper bound for the covering number
$M(\Ps_m(\Gamma), \sqrt{x}/24; \ell_f)$. The set $\Ps_m(\Gamma)$ is
much smaller than $S_m(P, \omega)$, however, and the following direct
argument gives the simpler upper bound:
%
\begin{equation}\label{diffco}\quad
M\bigl(\Ps_m(\Gamma), \sqrt{x}/24; \ell_f\bigr) \leq M\bigl(\Ps_m(\Gamma),
\sqrt{x}/24; \ell_H\bigr) \leq\biggl(1 +
\frac{48\Gamma}{\sqrt{x}} \biggr)^{md}.
\end{equation}
To see this, let $\mathfrak{C}$ be an $\epsilon:=
\sqrt{x}/24$-covering subset of $B_d(0, \Gamma)$ under the Euclidean
metric, $\euc_d$. Then $\mathfrak{H} := \{\operatorname{conv}(S) \dvtx S
\subseteq\mathfrak{C}\mbox{ with } |S| \leq m \}$ forms an
$\epsilon$-covering
subset of $\Ps_m(\Gamma)$ under $\ell_H$. Indeed, if $K =
\operatorname{conv}\{a_1,\ldots, a_m\}$ with $a_1,\ldots, a_m \in B_d(0,
\Gamma)$ is an arbitrary element set in $\Ps_m(\Gamma)$, then we can
choose $a_1',\ldots, a_m' \in\mathfrak{C}$ with $\euc_d(a_i,
a_i') \leq\epsilon$ for each $i$. It is then easy to see
[using~(\ref{althaus})] that the Hausdorff distance between $K$ and
$\operatorname{conv}\{a_1',\ldots, a_m'\}$ is at most $\epsilon$. It is
evident that the cardinality of $\mathfrak{H}$ is at most
$|\mathfrak{C}|^m$. A~standard volumetric argument shows that the
cardinality of $\mathfrak{C}$ can be taken to be smaller than
$(1+(2\Gamma/\epsilon))^d$ which proves~(\ref{diffco}).

Combining~(\ref{akka}) and~(\ref{diffco}), we get
\[
\varrho_x \leq4 \min\biggl(1, \biggl(1+\frac{48\Gamma}{\sqrt{x}}
\biggr)^{md} \exp\biggl(-\frac{nx}{2304\Gamma^2} \biggr) \biggr).
\]
The left-hand side of~(\ref{gyopoleq}) can now be bounded, for every
$\Lambda> 0$, by
\begin{eqnarray*}
\int_0^{\infty} \varrho_x \,dx &=& \int_0^{\Lambda} \varrho_x \,dx +
\int_{\Lambda}^{\infty} \varrho_x \,dx \\
&\leq& 4 \Lambda+ 4 \int_{\Lambda}^{\infty}
\biggl(1+\frac{48\Gamma}{\sqrt{x}} \biggr)^{md} \exp
\biggl(-\frac{nx}{2304\Gamma^2} \biggr) \,dx \\
&\leq& 4 \Lambda+ 9216 \frac{\Gamma^2}{n}
\biggl(1+\frac{48\Gamma}{\sqrt{\Lambda}} \biggr)^{md} \exp\biggl(-
\frac{n\Lambda}{2304 \Gamma^2} \biggr).
\end{eqnarray*}
The required bound~(\ref{gyopoleq}) is now easily deduced by using
the above inequality with
\[
\Lambda:= 2304 \frac{\Gamma^2}{n} md \log\bigl(1 + 48 \sqrt{n}
\bigr).
\]
\upqed
\end{pf}
\end{appendix}



\printaddresses

\end{document}